\documentclass[11pt]{amsart}

\usepackage{amssymb}
\usepackage{amsfonts,amsmath,amssymb,mathrsfs,verbatim}
\usepackage[dvips]{graphicx}
\usepackage{graphics}
\usepackage{psfrag}
\usepackage{hyperref}
\usepackage{marginnote}
\usepackage{xcolor}
\usepackage[dvips]{graphicx}
\usepackage{pgfplots}
\pgfplotsset{compat=1.16}
\usetikzlibrary{arrows}
\usepackage{tikz}
  \usetikzlibrary{shapes,backgrounds}
  \usetikzlibrary{calc}
  \usetikzlibrary{matrix}
\usepackage{caption}
\usepackage[margin=3cm]{geometry}
\usepackage{multicol}

\newcommand{\x}{\mathbf{x}}

\newtheorem{thm1}{Theorem}[section]
\newtheorem{lem1}[thm1]{Lemma}
\newtheorem{cor1}[thm1]{Corollary}

\newtheorem{prop1}[thm1]{Proposition}

\theoremstyle{definition}
\newtheorem{ex1}[thm1]{Example}
\newtheorem{def1}[thm1]{Definition}
\newtheorem{rem1}[thm1]{Remark}


\begin{document}

\title{Unimodular toric ideals of graphs}

\author[Ch. Tatakis]{Christos Tatakis}
\address{Department of Mathematics, University of Western Macedonia, 52100 Kastoria, Greece}
\email{chtatakis@uowm.gr}

\begin{abstract}
We give a necessary and sufficient graph-theoretic characterization of toric ideals of graphs that are unimodular. As a direct consequence, we provide the structure of unimodular graphs by proving that the incidence matrix of a graph $G$ is unimodular if and only if any two odd cycles of $G$ intersect.
\end{abstract}

\maketitle


\section{Introduction}

In the literature there are several results describing graphs that their toric ideals have a certain algebraic property, for instance, normality, complete intersection, robustness, generalized robustness, strongly robustness, generated by quadrics, quadratic Gr\"obner bases, Koszulness, see \cite{BOOCHER3, Nacho-T1, Gi-R-Vega, Gi-R-Vill, OH2, Ohsugi, Simis2, Tata1, Tata-Thom2, Vill2}. The goal of this paper is to classify the toric ideals of graphs that are unimodular.

For an integer matrix $A$, where $\textrm{rank}(A)=d$, the matrix $A$ is called unimodular if and only if all nonzero $d\times d$-minors of $A$ have the same absolute value. Unimodular matrices are of high interest for a lot of areas such as algebraic statistics, commutative algebra, algebraic geometry, integer programming, e.t.c., see for instance, \cite{Akbari, Bayer-Pop-Stu, Bern-ONeil, Bern-Sull, Camion, Sull2, Commoner, Conca, Sull3, Haase, OH3, Ohs-Tsu}. More precisely, in algebraic statistics for unimodular matrices it is easy to solve the integer programs that arise when evaluating whether individual entries of a data table are secure or when performing sequential importance sampling$;$ see \cite{Sull2, Sull}. In addition, the property of unimodularity is studied and completely characterized for hierarchical and binary hierarchical models, see \cite{Bern-ONeil, Bern-Sull}. Additionally, for unimodular matrices the Graver basis and their Markov bases are very easy to compute, see \cite{Sull3, ST}.

In algebraic geometry, it is known that if a matrix $A$ is unimodular, then the secondary polytope and the state polytope of $I_A$ coincide, as it holds for the Gr\"obner fan and the secondary fan of $I_A$, see [\cite{ST}, Proposition 8.15]. 

In commutative algebra we know the existence of a strong connection between square-free initial ideals of a toric ideal $I_A$ and unimodular regular triangulations of the edge polytope of $A$, see \cite{ST}. Sturmfels proved that a matrix is unimodular if and only if all initial ideals of its corresponding toric ideal are square-free, see [\cite{ST}, Remark 8.10]. It is also worth noting the connection between unimodular matrices and normality. In fact, for any unimodular matrix $A$, the corresponding semigroup ring $\mathbb{K}[A]$ is normal, while the converse is not true, see \cite{OH3}. A necessary condition for the corresponding toric ideal $I_A$ to have a square-free initial ideal is the normality of $\mathbb{K}[A]$, which was characterized combinatorially by Ohsugi and Hibi \cite{OH2} and Simis, Vasconcelos and Villarreal \cite{Simis2}. Sturmfels proved that if $I_A$ admits a square-free initial ideal with
respect to some term order, then $\mathbb{K}[A]$ is normal, see [\cite{ST}, Proposition 13.15], while by a well-known result of Hochster, see \cite{HOC}, we have that if $\mathbb{K}[A]$ is normal, then it is Cohen-Macaulay. Combining all the above, we conclude that for a unimodular toric ideal the corresponding polynomial ring is normal and thus Cohen-Macaulay$;$ a remark which increases the importance and interest of the classification of unimodular matrices.

In computer science and integer programming, there are a lot of applications through the totally unimodular matrices, a subclass of unimodular matrices. A matrix $A$ is totally unimodular if every square submatrix of $A$ has determinant $0$ or $\pm1$. It is known that the incidence matrix of a graph $G$ is totally unimodular if and only if $G$ is a bipartite graph. Totally unimodular matrices are very well behaved because they always define polytopes with integer vertices. Their significance and importance stem from a lot of applications such as min-cost perfect matching in bipartite graphs (assignment problem), maximum weight matching in bipartite graphs, flow problems e.t.c., see, for instance, \cite{Akbari, Camion, Commoner, Seymour}.

The main result of the present manuscript characterizes completely the graphs giving rise to unimodular toric ideals$;$ see Theorem \ref{maintheorem}. In terms of matrices, our main result is as follows.

\begin{thm1}\label{main-intro} The incidence matrix of a connected graph $G$ is unimodular if and only if any two odd cycles of $G$ intersect.
\end{thm1}

In order to prove the above result, we use the theory of toric ideals of graphs and their toric bases$;$ the set of circuits of a toric ideal, its universal Gr\"obner basis and its Graver basis. Finally, the above result leads us to give a structural way for the graphs whose toric ideals are unimodular, see Theorem \ref{corollary}. With the last result, we are able to describe all unimodular toric ideals of graphs, see Section \ref{last}.

\section{Background}

Let $A=\{\textbf{a}_1,\ldots,\textbf{a}_m\}\subseteq \mathbb{N}^n$
be a finite set of non-zero vectors  and
$\mathbb{N}A:=\{l_1\textbf{a}_1+\cdots+l_m\textbf{a}_m \ | \ l_i \in\mathbb{N}\}$ the 
corresponding affine semigroup. We grade the
polynomial ring $\mathbb{K}[x_1,\ldots,x_m]$ on an arbitrary field $\mathbb{K}$ by setting $\deg_{A}(x_i)=\textbf{a}_i$ for
$i=1,\ldots,m$. For $\textbf{u}=(u_1,\ldots,u_m) \in \mathbb{N}^m$,
we define the $A$-\emph{degree} of the monomial $\textbf{x}^{\textbf{u}}:=x_1^{u_1} \cdots x_m^{u_m}$
to be $\deg_{A}(\textbf{x}^{\textbf{u}}):=u_1\textbf{a}_1+\cdots+u_m\textbf{a}_m
\in \mathbb{N}A$. The
\emph{toric ideal} $I_{A}$ associated to $A$ is the prime ideal generated by all the $A$-homogeneous binomials, i.e., $$I_A=\langle \textbf{x}^{\textbf{u}}- \textbf{x}^{\textbf{v}}\ \textrm{such that}\ \deg_{A}(\textbf{x}^{\textbf{u}})=\deg_{A}(\textbf{x}^{\textbf{v}})\rangle.$$ 

Some of the very important toric bases of a toric ideal are its Graver basis, its universal Gr\"obner basis and the set of the circuits of the ideal. A binomial $\textbf{x}^{\textbf{u}}-
\textbf{x}^{\textbf{v}}$ in $I_A$ is called {\it primitive} if
there is no other binomial
 $\textbf{x}^{\textbf{w}}- \textbf{x}^{\textbf{z}}$ in $I_A$,
such that $\textbf{x}^{\textbf{w}}$ divides $
\textbf{x}^{\textbf{u}}$ and $\textbf{x}^{\textbf{z}}$ divides $
\textbf{x}^{\textbf{v}}$. The set of primitive binomials, which is finite, is the
{\it Graver basis} of $I_A$ and is denoted by $Gr_A$.   
 The {\it universal Gr\"{o}bner basis}
 of an ideal $I_A$, is denoted by $\mathcal{U}_A$ and is defined as the union of all reduced Gr\"obner bases $G_\prec$ of $I_A$, as $\prec$ runs over all term orders. It is a finite subset of binomials in $I_A$ and is a Gr\"obner basis for the ideal with respect to all term orders, see \cite{ST}. The support of a monomial ${\bf{x^u}}$ of $\mathbb{K}[x_1,\ldots,x_m]$
is $\textrm{supp}({\bf{x^u}}):=\{i\ | \ x_{i}\ \mbox{ divides}\ {\bf{x^u}}\}$ and
the support of a binomial $B={\bf{x^u}}-{\bf{x^v}}$ is
$\textrm{supp}(B):=\textrm{supp}({\bf{x^u}})\cup \textrm{supp}({\bf{x^v}})$.
An irreducible non-zero binomial is called a circuit if it has minimal support. Equivalently, in terms of matrices, for an integer matrix $A$, a non-zero element ${\bf {u}}\in\ker_{\mathbb{Z}}A$ is called a circuit of $A$ if its non-zero entries are relatively prime and there is no other non-zero element  ${\bf {v}}\in\ker_{\mathbb{Z}}A$ such that $\textrm{supp}({\bf {v}})\subseteq \textrm{supp}({\bf {u}})$. The set of circuits of a toric ideal $I_A$ is denoted by $\mathcal{ C}_A$.

The relation between the above toric bases was given by B. Sturmfels.

\begin{prop1}\label{CUG}\cite[Proposition 4.11]{ST} For any toric ideal $I_A$ it holds:
$$\mathcal{ C}_A\subseteq  \mathcal{ U}_A\subseteq Gr_A$$
\end{prop1} 

For a deeper treatment of toric bases, see \cite{Hi-Ohs, RTT, ST, Tata-Thom1, Vill1}.

In the next chapters,  $G$ is a simple, connected, undirected, and finite graph, for which we denote by $V(G)$ the set of its vertices, and let $E(G)=\{e_{1},\ldots,e_{m}\}$ be
the set of its edges. Let $\mathbb{K}[e_{1},\ldots,e_{m}]$ be the polynomial ring in the $m$ variables $e_{1},\ldots,e_{m}$ on an arbitrary field $\mathbb{K}$.  We
will associate each edge $e=\{v_{i},v_{j}\}\in E(G)$  with the element
$a_{e}=v_{i}+v_{j}$ in the free abelian group $ \mathbb{Z}^n $, with the basis the set of vertices
of $G$, where $v_{i}=(0,\ldots,0,1,0,\ldots,0)$ be the vector with 1 in the $i-$th coordinate of $v_{i}$. By $I_{G}$ we denote the toric ideal $I_{A_{G}}$ in $\mathbb{K}[e_{1},\ldots,e_{m}]$, where $A_{G}=\big\{a_{e}\ | \ e\in E(G)\big\}\subseteq \mathbb{Z}^n $.

In order to better describe the toric ideal of graph and its toric bases, we need some basic elements of graph theory. A {\it walk} connecting $u \in V(G)$ and
$u' \in V(G)$ is a finite sequence of vertices of graph $w=(u = u_0, u_1, \ldots, u_{\ell-1}, u_\ell = u')$,
with each $e_{i_j}=\{u_{j-1},u_{j}\}\in E(G)$, for $j=1,\ldots,\ell$. 
The {\it length} of the walk $w$ is the number $\ell$ of its edges. An
{\it even} (respectively, {\it odd}) {\it walk} is a walk of even (respectively, odd) length. A walk
$w=(u_0,u_1\ldots,u_{\ell-1},u_{\ell})$
is called {\it closed} if $u_{0}=u_{\ell}$. A {\it cycle}
is a closed walk
$(u_0,u_{1},\ldots,u_{\ell-1},u_{\ell})$ with
$u_{k}\neq u_{j},$ for every $ 1\leq k < j \leq \ell$, while a {\it path} is a walk of the graph where all its vertices are distinct. A chord of a walk $w$ is an edge of the graph $G$ that joins two non-adjacent vertices of the walk $w$. A walk $w$ is called chordless if it does not have chords. Finally, a {\em cut edge} (respectively, {\em cut vertex}) is an edge (respectively, vertex) of
the graph whose removal increases the number of connected
components of the remaining subgraph.  A graph is called {\em
biconnected} if it is connected and does not contain a cut
vertex. A {\em block} is a maximal biconnected subgraph of a given
graph $G$.

Consider an even closed walk $w=(u_{0},u_{1},u_{2},\ldots,u_{2s-1},u_{2s} = u_0)$ of length $2s$ with $e_{i_j}=\{u_{j-1},u_{j}\}\in E(G)$, for $j=1,\ldots,2s$. The binomial $B_w = e_{i_1}e_{i_3} \cdots e_{i_{2s-1}}  - e_{i_2} e_{i_4} \cdots e_{i_{2s}}$
belongs to the toric ideal $I_G$. In fact, Villarreal proved that $$I_G = \langle B_w  \ \vert \ w \text{ is an even closed walk}\rangle, $$ 
that is, the toric ideal $I_G$ is generated by the binomials corresponding to even closed walks of the graph $G$, see \cite{Vill1, Vill3}.

In the case of toric ideals of graphs, all the toric bases are known, see \cite{Hi-Ohs, RTT, Tata-Thom1, Vill1}. The following theorems determine the form of the circuits and the
primitive binomials of a toric ideal of a graph $G$. For the sake of brevity, we refer the reader to the corresponding articles. Villarreal gave a necessary and sufficient
characterization of the circuits (that is, the set $\mathcal{C}_G$). For convenience by $\bf{w}$
we denote the subgraph of $G$ with vertices the vertices of the
walk and edges the edges of the walk $w$. Note that $\bf{w}$
is a connected subgraph of $G$.
\begin{thm1}\cite[Proposition 4.2]{Vill1}\label{circuit} Let $G$ be a graph and let $W$ be a connected subgraph of $G$.
The subgraph $W$ is the graph  ${\bf w}$ of a walk $w$ such that  $B_w$ is a circuit
 if and only if
\begin{enumerate}
  \item[($c_1$)] $W$ is an even cycle or
  \item[($c_2$)] $W$ consists of two odd cycles intersecting in exactly one vertex or
  \item[($c_3$)] $W$ consists of two vertex-disjoint odd cycles joined by a path.
\end{enumerate}
\end{thm1}
 
 From \cite{Hi-Ohs} we also know the form of the primitive walks of a graph $G$.

\begin{lem1}\label{primitOH}\cite[Lemma 3.2]{Hi-Ohs} If $B_w$ is primitive, then $w$ has one of the
following forms:
\begin{enumerate}
  \item[($p_1$)] $w$ is an even cycle or
  \item[($p_2$)] $w$ consists of two odd cycles intersecting in exactly one vertex or
  \item[($p_3$)] $w=(c_1,w_1,c_2,w_2)$ where $c_1,c_2$ are odd vertex disjoint cycles and $w_1,w_2$
  are walks which combine a vertex $v_1$ of $c_1$ and a vertex $v_2$ of $c_2$.
\end{enumerate}
\end{lem1}

In the following example, we illustrate the similarities and differences between circuits and primitive elements of the toric ideals of graphs, the understanding of which plays a crucial role in the next chapter.

\begin{ex1}{\rm

By Sturmfels, we know that $\mathcal{C}_G\subseteq Gr_G$. The converse inclusion also holds in the case that the walk $w$ has either the form $(p_1)$ or the form $(p_2)$ of Lemma \ref{primitOH}, see Figure \ref{ci=primi}. 

\begin{figure}[h]
\begin{center}
\includegraphics[scale=1.2]{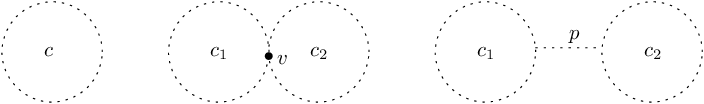}
\caption{The cases that circuits and primitive elements coincide$;$ $c$ is an even cycle, $v$ a vertex, $c_1,c_2$ are odd cycles and $p$ a path}
\label{ci=primi}
\end{center}
\end{figure}

The figure on the left hand consists of an even cycle $c$, the figure on the middle consists of two odd cycles intersecting in exactly one vertex $v$ of $G$, and the last figure on the right hand consists of two odd disjoint cycles $c_1,c_2$ joined by a path $p$ of length at least one. From Theorem \ref{circuit} and Lemma \ref{primitOH} the corresponding binomials are circuits and elements of the Graver basis.

However, when the walk $w$ has the form $(p_3)$ of Lemma \ref{primitOH}, the corresponding binomial $B_w$ instead of being sometimes primitive, it is not a circuit, see Figure \ref{priminotcir}.

\begin{figure}[h]
\begin{center}
\includegraphics[scale=1.2]{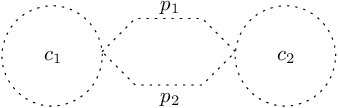}
\caption{The case that a primitive element is not a circuit$;$ $c_1,c_2$ are odd cycles which are joined by two disjoint paths $p_1$ and $p_2$ of the same parity$;$ i.e. they are both even or odd}
\label{priminotcir}
\end{center}
\end{figure}

}
\end{ex1}

The next theorem by  E.~Reyes et all, describes  the form of the underlying graph of a primitive walk and thus gives us the Graver basis $Gr_G$ of the ideal $I_G$. In the following theorem we can verify (see condition (2b)), that the binomial $B_w$ where $w$ is the walk which corresponds to the walk of Figure \ref{priminotcir} belongs to the Graver basis of the corresponding toric ideal.

\begin{thm1}\cite[Corollary 3.3]{RTT} \label{primitive-graph}
Let $G$ be a graph, and let $W$ be a connected subgraph of $G$.
The subgraph $W$ is the graph  ${\bf w}$ of a primitive walk $w$
 if and only if \begin{enumerate}
  \item  $W$ is an even cycle or
  \item  $W$ is not biconnected and
\begin{enumerate}
  \item every block of $W$ is a cycle or a cut edge and
  \item every cut vertex of $W$ belongs to exactly two blocks and separates the graph in two parts, 
  the total number of edges
of the  blocks that are cycles in each part is odd.
\end{enumerate}
\end{enumerate}
\end{thm1}

\section{unimodular toric ideals of graphs}

We start this section by setting how the set of the circuits and the Graver basis of a toric ideal behave with respect to elimination of variables.
 
\begin{prop1}\cite[Proposition 4.13]{ST}\label{elimination} 
Let $A$ be a finite set of positive integers. If $A' \subseteq A$ is not empty, then
\begin{itemize}
\item[($\alpha$)] $\mathcal C_{A'} = \mathcal C_A \cap \mathbb K[\x_{A'}]$,
\item[($\beta$)] $\operatorname{Gr}_{A'} = \operatorname{Gr}_A \cap \mathbb K[\x_{A'}]$.
\end{itemize}
where $\mathbb K[\x_{A'}] :=  \mathbb K[x_i \mid a_i \in A']$.
\end{prop1}

The main goal of this manuscript is to describe in graph-theoretical terms the toric ideals of graphs that are unimodular. Unimodularity is a strong property that an integral matrix $A$ could satisfy. We could define a unimodular matrix if the entries in each circuit of $A$ (and also the elements of its Graver basis) have entries either $0$ or $\pm 1$. The most common definition is the following$;$ see \cite{Bern-Sull}.

\begin{def1} If $rank(A)=d$, the matrix $A$ is called unimodular if and only if all non-zero $d\times d$-minors of $A$ have the same absolute value.
\end{def1}

We say that a toric ideal $I_A$ is unimodular if the corresponding matrix $A$ is unimodular. The notion of square-free ideals plays a key role for unimodular matrices. We recall that a monomial ${\bf{x^u}}$ is square-free if every coordinate of ${\bf{u}}$ is 0 or 1. A binomial is square-free if its monomials are square-free. An ideal is square-free if its generators are square-free. 

\begin{def1} A graph $G$ is called unimodular if its incidence matrix (i.e. its corresponding toric ideal) is unimodular.
\end{def1}

 In order to examine the property of unimodularity for the toric ideals of graphs, we set the following properties.

\begin{thm1}\cite[Remark 8.10]{ST}\label{iffAuni} A matrix $A$ is unimodular if and only if all initial ideals of the toric ideal $I_A$ are square-free.
\end{thm1}

Unimodular matrices also have the following important property.

\begin{prop1}\cite[Proposition 8.11]{ST}\label{Sturm-prop} Let $A$ be a unimodular matrix and let $I_A$ be its corresponding toric ideal. The set of the circuits $C_A$ equals its Graver basis $Gr_A$.
\end{prop1}

We remark that for the converse statement of Proposition \ref{Sturm-prop}, we need the circuits of the toric ideal to be square-free, as we prove in the next proposition.

\begin{prop1}\label{converse-uni} A matrix $A$ is unimodular if and only if the set of the circuits equals the Graver basis and they are square-free.
\end{prop1}
\begin{proof}
If the matrix $A$ is unimodular the result follows by Proposition \ref{Sturm-prop} and the fact that, by definition, each circuit of $A$ has entries either $0$ or $\pm 1$. 

For the converse statement, by hypothesis and Proposition \ref{CUG} it follows that $\mathcal{ C}_A=  \mathcal{ U}_A=Gr_A$. Thus, the binomials of the universal Gr\"obner basis of $I_A$ are square-free, which means that all its initial ideals of $I_A$ are square-free. The result follows from Theorem \ref{iffAuni}.
\end{proof}

From the previous proposition, we can prove that the unimodularity between toric ideals is a hereditary property, that is, it is closed when taking subsets$;$ see also \cite{Bern-Sull}.

\begin{prop1}\label{hereditary} Let $B\subseteq A$. If $I_A$ is a unimodular toric ideal then $I_B$ is a unimodular toric ideal.
\end{prop1}
\begin{proof}
Let $B\subseteq A$. From Proposition \ref{converse-uni} we have to prove that $\mathcal{C}_B=Gr_B$ and all binomials of $\mathcal{C}_B$ are square-free.

From Proposition \ref{CUG} we have $C_B\subseteq Gr_B$. For the converse inclusion, let $f\in Gr_B$ then $\textrm{supp}(f)\subseteq B$, and from Proposition \ref{elimination} we have $f\in Gr_A$. Since $I_A$ is unimodular, it follows that $Gr_A=\mathcal{C}_A$ and thus $f\in \mathcal{C}_A$. It follows that the binomial $f$ has minimal support in the set $A$ and therefore has minimal support in any subset of $A$ that contains $\textrm{supp}(f)$, as the set $B$, that is, $\textrm{supp}(f)\subseteq B\subseteq A$ and thus $f\in\mathcal{C}_B$. It follows that $\mathcal{C}_B=Gr_B$.

Also, since $I_A$ is unimodular, it follows that the binomials of its Graver basis (which are also circuits) are square-free. From Proposition \ref{elimination} it follows that the binomials of $Gr_B$ (and of $\mathcal{C}_B$) are also square-free. The result follows.
\end{proof}

An immediate application of the above proposition is the following useful result.

\begin{cor1}Let $G$ be a graph. The ideal $I_G$ is unimodular if and only if for every connected component $H$ of $G$ the ideal $I_H$ is unimodular.
\end{cor1}
\begin{proof}
If the toric ideal $I_G$ is unimodular, the result follows from Proposition \ref{hereditary}.

For the converse statement, by the definition of a toric ideal of a graph $G$, every generator of the ideal and thus every binomial of its Graver basis belong to a connected component of the graph $G$. The result follows.
\end{proof}

The above result allows us to examine the problem of unimodular graphs for the case of connected graphs. By Proposition \ref{hereditary}, we have the following corollary.

\begin{cor1}\label{blocks} Let $G$ be a connected graph. If the ideal $I_G$ is unimodular then for every block $H$ of $G$ the ideal $I_H$ is unimodular.
\end{cor1}

The converse of the above corollary is not true as we can see in the next example.

\begin{ex1}\label{example}{\rm

In this example, we see that the converse statement of Corollary \ref{blocks} does not hold. The graph in Figure \ref{notunimodular} consists of two non bipartite blocks$;$ $$B_1=\{x_1,x_2,x_3,x_4,x_5\} \ \textrm{and}\ B_2=\{x_6,x_7,x_8,x_9,x_{10}\}.$$

\begin{figure}[h]
\begin{center}
\includegraphics[scale=1.2]{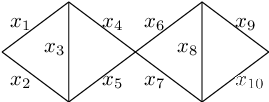}
\caption{A not unimodular graph. All its blocks are unimodular.}
\label{notunimodular}
\end{center}
\end{figure}

It is easy to check that

$$\mathcal{C}_{B_1}=Gr_{B_1} =  \langle   x_{1}x_{5}-x_{2}x_{4} \rangle \ \textrm{and}\ \mathcal{C}_{B_2}= Gr_{B_2} =  \langle   x_{6}x_{10}-x_{7}x_{9} \rangle$$

As we see all the binomials are square-free. By Proposition \ref{converse-uni}, it follows that both $B_1,B_2$ are unimodular graphs.

However, the whole graph $G=B_1\cup B_2$ is not unimodular. By computations we have that

\begin{eqnarray}
\mathcal{C}_G &=&  Gr_G  =  \langle  x_3x_2x_6^2x_{10}-x_1x_5^2x_9x_8, x_4x_5x_8-x_3x_6x_7, x_4x_5x_9x_8-x_3x_6^2x_{10}, \nonumber\\
& &  x_4x_5x_{10}x_8-x_3x_7^2x_{9}, x_3x_1x_7^2x_{9}-x_2x_4^2x_{10}x_8, x_{1}x_{5}-x_{2}x_{4}, x_{6}x_{10}-x_{7}x_{9},  \nonumber\\
& & x_3x_2x_7x_6-x_1x_5^2x_8, x_3x_1x_7x_6-x_2x_4^2x_8, x_3x_1x_6^2x_{10}-x_2x_4^2x_9x_8, x_3x_2x_7^2x_{9}-x_1x_5^2x_{10}x_8 \rangle \nonumber
\end{eqnarray}

We remark the existence of the binomials $$ x_3x_2x_6^2x_{10}-x_1x_5^2x_9x_8, x_3x_1x_7^2x_{9}-x_2x_4^2x_{10}x_8, x_3x_1x_6^2x_{10}-x_2x_4^2x_9x_8, x_3x_2x_7^2x_{9}-x_1x_5^2x_{10}x_8.$$ Any of the above binomials belong also to the universal Gr\"obner basis of the ideal, which give us not square-free binomials in an initial ideal of $I_G$.  By Proposition \ref{iffAuni} it follows that the ideal is not unimodular. Note that someone can conclude the non-unimodularity by applying Proposition \ref{converse-uni}.

}
\end{ex1}

We note that instead of a small graph (as in the previous example), the computations are complicated$;$ for more complicated graphs, the corresponding computations become extremely difficult. Our aim is to give a structural way for unimodular graphs to avoid all the corresponding difficulties. Next, we recall a useful definition for our theorem.

\begin{def1} Let $G$ be a graph. We say that $G$ has the strong odd cycle property if any two odd cycles intersect.
\end{def1}

For example, the complete graph on the $n$ vertices $K_n$, has the strong odd cycle property for $n\leq 5$ but it has not for any $n>5$, because of the existence of two triangles not intersecting.

Next we state the main result of this manuscript which characterizes completely when a toric ideal of a graph $G$ is unimodular.

\begin{thm1}\label{maintheorem} Let $G$ be a connected graph. The toric ideal $I_G$ is unimodular if and only if $G$ has the strong odd cycle property.
\end{thm1}

Equivalently, as we stated in the introduction, in graph theory terms, the above theorem can be written as the incidence matrix of a graph $G$ is unimodular if and only if any two odd cycles of $G$ intersect$;$ see Theorem \ref{main-intro}. For example, it follows that all the connected graphs of four vertices are unimodular.

In order to prove Theorem \ref{maintheorem}, it is enough to prove the following equivalent theorem. The following result is a different approach to Theorem \ref{maintheorem}, which gives us a better view of the unimodular graphs in the case that the graphs have either two or more non bipartite blocks. As we prove below, in the above case, there exists a common vertex for all odd cycles of $G$ through which they are passing. This vertex (which is called a link vertex of the graph) is a cut vertex of the graph $G$. A vertex $v$ of $G$ is called a link vertex if every odd cycle of $G$ passes through $v$. Note that due to Corollary \ref{blocks} the following theorem holds for each connecting component of a graph $G$.

\begin{thm1}\label{corollary}
Let $G$ be a connected graph. The toric ideal $I_G$ is unimodular if and only if exactly one of the following holds:
\begin{itemize}
\item[($\alpha$)] All blocks of $G$ are bipartite.
\item[($\beta$)] All blocks of $G$ are bipartite except one that has the strong odd cycle property.
\item[($\gamma$)] All blocks of $G$ are bipartite except $s\geq 2$ blocks. In this case $G$ has a link vertex $x$.
\end{itemize}
\end{thm1}
\begin{proof} 

$(\Longleftarrow )$ By hypothesis, the graph $G$ either is bipartite (case ($\alpha$)) or any two odd cycles intersect (cases $(\beta), (\gamma)$). Let $B_w$ be an element of the Graver basis of $I_G$. Combining Lemma \ref{primitOH} and the form of $G$, we have that $w$ is either an even cycle of $G$ or consists of two odd cycles $c_1,c_2$ such that $V(c_1)\cap V(c_2)=\{v\}$, where $v$ is a vertex of $G$. In any case, according to Theorem \ref{circuit}, $B_w$ is also a circuit, which means that $Gr_G\subseteq C_G$ and therefore by Proposition \ref{CUG} we have $C_G=Gr_G$. 

We claim that $B_w$ is square-free. Suppose not. Since $C_G=Gr_G$, by Theorem \ref{circuit} it follows that $w$ consists of two disjoint odd cycles joined by a path of length at least one (since the graph is also connected). A contradiction arises due to the fact that the graph has the strong odd cycle condition. By Proposition \ref{converse-uni}, it follows that $I_G$ is unimodular.

\textit{$(\Longrightarrow )$} Let $I_G$ be a unimodular toric ideal and let $s$ be the number of non bipartite blocks of $G$. If $s=0$ the result follows. Suppose now that the graph $G$ is not bipartite, that is, $s\geq 1$. First, we will prove that $G$ has the strong odd cycle property.

We suppose that there exist at least two odd disjoint cycles of $G$ and let them be $c_1=(e_{1,1},\ldots,e_{1,2k+1})$ and $c_2=(e_{2,1},\ldots,e_{2,2m+1})$, where $V(c_1)\cap V(c_2)=\emptyset$. Since the graph is connected, there exists a path $p=(\epsilon_1,\ldots,\epsilon_n)$ with $n\geq1$ that connects $c_1$ and $c_2$. We consider the walk $w=(c_1,p,c_2,-p)$. The corresponding binomial has the form $$B_w=\epsilon_2^2\cdots\epsilon_{n}^2\prod_{i=0}^{i=k}e_{1,2i+1}\prod_{j=0}^{i=m}e_{2,2j+1}-\epsilon_1^2\cdots\epsilon_{n-1}^2\prod_{i=1}^{j=k}e_{1,2i}\prod_{j=1}^{j=m}e_{2,2j}$$ 

By Theorem \ref{circuit} $(c_3)$, $B_w$ is a circuit (and thus an element of the universal Gr\"obner basis of $I_A$) which is not square-free, contradicting the fact that the ideal is unimodular.

Suppose now that the graph $G$ has $s\geq 2$ non bipartite blocks. It remains to prove that $G$ has a link vertex.

Let $B_1$ and $B_2$ be two different non bipartite blocks of $G$ and let $c_1=(v_1,\ldots,v_{2k+1})$ and $c_2=(u_1,\ldots,u_{2l+1})$ be correspondingly two odd cycles of these blocks. By our previous claim, the graph $G$ has the strong odd cycle property, thus we have $V(c_1)\cap V(c_2)\neq \emptyset$. Since the cycles belong to different blocks, it follows that $V(c_1)\cap V(c_2)=\{v\}$, where $v\in V(G)$ and without loss of generality we suppose that $v=v_{1}=u_1$. We aim for the vertex $v$ to be a link vertex of $G$. 

Suppose not. Then there exists an odd cycle $c_3$ of $G$ such that $v\notin V(c_3)$. Since the graph $G$ has the strong odd cycle property, we have $V(c_1)\cap V(c_3)\neq \emptyset$ and $V(c_2)\cap V(c_3)\neq \emptyset$. Let $i,j$ be the smallest possible values such that $v_i, u_j\in V(c_3)$, where $v_i,u_j$ are different from the vertex $v$. It follows that there exist at least two disjoint paths $p_1, p_2$ of the graph $G$ which join the vertices $v_i, u_j;$

$$p_1=(v_i,v_{i-1},v_{i-2},\ldots,v_{1}=v=u_1,u_2,\ldots,u_{j-1},u_j)$$ and the path $p_2$ which consists of vertices of the cycle $c_3$.

The contradiction arises because the vertices $v_i$ and $u_j$ belong to different blocks of $G$.

\end{proof}

\begin{ex1}{\rm 

In this example we would like to present the differences between the unimodular graphs and how they look like in the cases of either they have one or with more than one non bipartite blocks, as we mention them in Theorem \ref{corollary}, see Figure \ref{unimodular}. From Theorem \ref{maintheorem} it follows that both toric ideals are unimodular. Both figures are non bipartite graphs which have the strong odd cycle property. Note that in the figure on the right, the graph has four non bipartite blocks and the vertex $v$ is a link vertex, which is a cut vertex of $G$.

\begin{figure}[h]
\begin{center}
\includegraphics[scale=1.2]{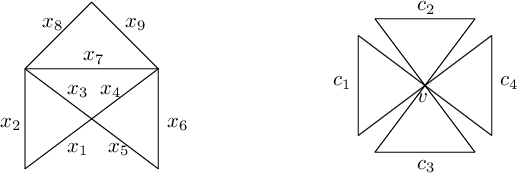}
\caption{Unimodular graphs. The graph on the left consists of one non bipartite block with the strong odd cycle property. The graph on the right consists of four odd cycles with one link vertex $v$.}
\label{unimodular}
\end{center}
\end{figure}

Differently, by computations for the toric ideal $I_G$, where $G$ is presented in Figure \ref{unimodular} on the left hand, we check that 
\begin{eqnarray}
\mathcal{C}_G\ =\ Gr_G \ =  &\langle &x_3x_9-x_4x_8, x_1x_7-x_2x_4, x_3x_6-x_5x_7, x_1x_6x_8-x_2x_5x_9, \nonumber\\
& & x_3x_2x_9-x_1x_8x_7, x_4x_6x_8-x_5x_7x_9, x_1x_3x_6-x_2x_4x_5 \ \rangle. \nonumber
\end{eqnarray}

while for the corresponding toric ideal of the graph on the right hand we check that 

$$\mathcal{C}_G=Gr_G=\langle B_w,\ w=(c_i,c_j),\ \textrm{where}\ i\neq j\in\{1,2,3,4\}\rangle$$

For both ideals, all the elements are square-free and by Proposition \ref{converse-uni} it follows that the ideals are unimodular.

}
\end{ex1}

\begin{ex1}{\rm 

We return to the graph of Example \ref{example}, see Figure \ref{notunimodular}. It is easy to check that it does not have the strong odd cycle property, since there are two disjoint odd cycles, that is, $c_1=(x_1,x_2,x_3)$ and $c_2=(x_{8},x_{9},x_{10})$. It follows from Theorem \ref{maintheorem} that the corresponding toric ideal is not unimodular. On the other hand, one can check that the graph has more that one non bipartite blocks that do not have a link vertex. Applying Theorem \ref{corollary}, we conclude the non unimodularity.

}
\end{ex1}

\section{The structure of unimodular graphs}\label{last}

Theorem \ref{corollary} is a structural algorithmic result instead of Theorem \ref{maintheorem}. The main advantage of Theorem \ref{corollary} is that we have a complete picture of unimodular graphs. In this way, we are able to construct as many (all) unimodular toric ideals of graphs as we want. 

For the construction of unimodular graphs, we recall that given a graph $H$, we call a path an $H$-path if it is nontrivial and meets $H$ exactly at its ends. If two paths have both an even or odd length, we say that they are of the same parity.

From our main result, we know that for all bipartite graphs the corresponding toric ideals are unimodular. In this case, it follows the above construction.

\begin{thm1} Let $G$ be a connected bipartite graph. $I_G$ is unimodular if and only if $G$ can be constructed from an even chordless cycle $G_0$, by successively adding $s$ $G_i$-paths by starting with $G_0$ and ending with $G_{s}=G$, where $i=0,\ldots,s-1$. The $G_i$-paths must be of the same parity with the parity of the path which joins their ends in the graph $G_i$.
\end{thm1}
\begin{proof}
By construction, the graphs $G_i$ are bipartite for all $i=0,\ldots,s$. The result follows from Theorem \ref{corollary}.
\end{proof}

For the non-bipartite case, the situation is much more complicated. The difficulties stem from the fact that for a toric ideal of a graph $G$, its minimal generators, and thus the elements of the toric bases of $I_G$ are much more complicated$;$ for more see \cite{RTT}. The advantage of Theorem \ref{corollary} is that it completely clarifies the situation. 

In order to describe, in graph-theoretical terms, the family of unimodular graphs with two or more non bipartite blocks, we need to introduce the notion of a flower-graph. 

\begin{def1} A graph $G$ is called a flower-graph if it consists of two or more odd chordless cycles $c_1,c_2,\ldots,c_k$ such that $V(c_1)\cap V(c_2)\cap \ldots V(c_k)=\{v\}$, where $v$ is a cut vertex of $G$. The vertex $v$ is called carpel, see the graph on the right of Figure \ref{unimodular}.
\end{def1}

By definition, it follows that the carpel $v$ is a link vertex of a flower-graph. Note that for any graph $G$ that has a link vertex, the graph $G\setminus \{v\}$ is bipartite$;$ thus it also holds for the flower-graphs. By Theorem \ref{corollary} we have the following result, which completes the construction of the unimodular toric ideal of graphs with at least two non bipartite blocks.

\begin{figure}[h]
\begin{center}
\includegraphics[scale=1.2]{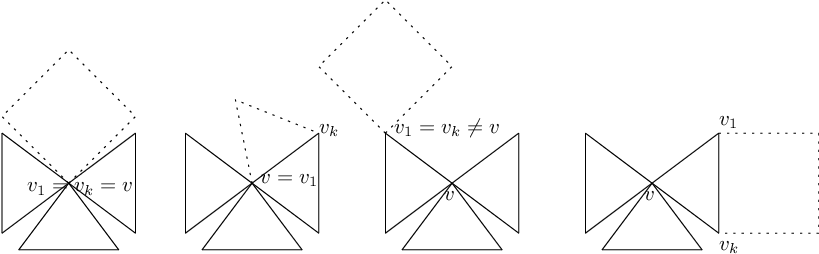}
\caption{The $H$-paths of Theorem \ref{construction!}. The first two figures correspond to the case ($\alpha$), the third figure correspond to the case ($\beta$), and the last one correspond to the case ($\gamma$)}.
\label{H-paths}
\end{center}
\end{figure}

\begin{thm1}\label{construction!} Let $G$ be a graph with two or more non bipartite biconnected blocks. $I_G$ is unimodular if and only if $G$ can be constructed from a flower-graph with carpel $v$, by successively adding $s$ $G_i$-paths, where $i=0,\ldots,s-1$, by starting with a flower-graph $G_0$ and ending with $G_{s}=G$, where the ends $v_1,v_k$ of each $G_i$-path belong to the same block of $G_{i}$, and the addition is as follows, see Figure \ref{H-paths}.
\begin{itemize}
\item[($\alpha)$] If $v_1=v$ or $v_k=v$, the $G_i$-path can be of any length,
\item[($\beta)$] if $v_1=v_k\neq v$, the $G_i$-path must be of even length,
\item[($\gamma)$] if $v_1,v_k$ are distinct and different from $v$, the $G_i$-path must be of the same parity as the parity of the path that joins $v_1$ and $v_k$ in the graph $G\setminus \{v\}$.
\end{itemize}

\end{thm1}

\begin{proof} Let $G$ be a graph with $\lambda\geq 2$ non bipartite biconnected blocks and let them be $B_1,\ldots,B_{\lambda}$.

$(\Longrightarrow )$ Let $c_1,\ldots,c_{\lambda}$ be odd chordless cycles of the blocks $B_1,\ldots,B_{\lambda}$ correspondingly. Since $I_G$ is unimodular, by Theorem \ref{corollary} we have that the graph $G$ has a link vertex, and let it be $v$. Due to the fact that the cycles $c_1,\ldots,c_{\lambda}$ belong to different blocks, we have that $V(c_1)\cap \ldots \cap V(c_{\lambda})=\{v\}$. It follows that the graph $G_0$ with edges $E(G_0)=E(c_1)\cup\ldots \cup E(c_{\lambda})$ and vertices $V(G_0)=V(c_1)\cup\ldots \cup V(c_{\lambda})$ is a flower-graph with carpel $v$ and it is a subgraph of $G$. 

If $G_0=G$ the result follows, and we suppose that it holds for any $\kappa\ < s$.  Let $p=(v_1,v_2,\ldots,v_k)$ be a $G_{s-1}$-path that joins two vertices $v_1$ and $v_k$ of $G_{s-1}$. There are two cases$;$ ($i$) $v_1=v$ or $v_k=v$ (see the first two figures on the left of Figure \ref{H-paths}) and ($ii$) the vertices $v_1,v_k$ are different from $v$ (see the last two figures on the right of Figure \ref{H-paths}). 

If the case ($i$) holds, since $v$ is a vertex of the $G_{s-1}$-path, any new cycle (and therefore for any new odd cycle) that is added in the graph $G_{s-1}$, passes through $v$. It follows that the vertex $v$ is a link vertex of the graph $G_s=G$ for any length (even or odd) of the $G_{s-1}$-path and the case ($\alpha$) follows. We consider now that the case ($ii$) holds. There are two cases$;$ either $v_1=v_k$ or $v_1\neq v_k$. If $v_1=v_k$, we suppose that the $G_{s-1}$-path has odd length. It follows that there exists the odd cycle $c=(v_1,v_2,\ldots,v_{k-1},v_k=v_1)$ that does not pass through the link vertex $v$, a contradiction arises due to the fact that the $I_G$ is unimodular and Theorem \ref{corollary}. Therefore the $G_{s-1}$-path has even length and the case ($\beta$) follows. For the last case, we suppose that $v_1\neq v_k$ and both of them differ from the link vertex $v$. Consider the bipartite graph $G_{s-1}\setminus \{v\}$, where $v_1,v_k\in V(G_{s-1}\setminus \{v\})$. If the parity of the $G_{s-1}$-path (let it be $p'$) differs from the parity of the path $p$ that joins the vertices $v_1,v_k$ in the graph $G_{s-1}\setminus \{v\}$, it follows that there exists an odd cycle $c'=(p',p)$ which does not pass though the link vertex $v$ in the graph $G_s=G;$ a contradiction arises similar to the previous case. Therefore the $G_{s-1}$-path has the same parity as $p$ and the case ($\gamma$) follows.

$(\Longleftarrow )$ Suppose that $G$ is constructed from a flower-graph with carpel $v$, by successively adding $s$ $G_i$-paths. By Theorem \ref{corollary}, it is enough to prove that $G$ has a link vertex. We will prove inductively on the number of $s$ $G_i$-paths that we added, that the carpel $v$ is a link vertex of $G$. 

For $s=0$, by definition of a flower-graph it follows that the carpel $v$ is a link vertex of $G_0=G$. Suppose that $v$ is a link vertex of the graph $G_{s-1}$, and we will prove that $v$ is a link vertex of the graph $G_s=G$.

By construction the graph $G_s$ arises from the graph $G_{s-1}$ by adding a $G_{s-1}$-path  of type of the cases that are described in $(\alpha), (\beta), (\gamma)$. If the $G_{s-1}$-path is of type that is described in the case ($\alpha$), we have that the vertex $v$ is a vertex of the $G_{s-1}$-path and therefore any new odd cycle of $G_s$ passes through $v$. Thus, the vertex $v$ is a link vertex of $G_s=G$. If we are in the case ($\beta$), obviously the $G_{s-1}$-path is an even cycle with one common vertex (the vertex $v_1=v_k$) with the graph $G_{s-1}$ and the vertex $v$ is a link vertex of $G_s$. For the case ($\gamma$), since the graph $G_{s-1}\setminus \{v\}$ is bipartite, and the $G_{s-1}$-path has the same parity with the path that joins $v_1$ and $v_k$ in the graph $G_{s-1}\setminus \{v\}$, it follows that the graph $G_s\setminus \{v\}$ is bipartite. Therefore the vertex $v$ is a link vertex of the graph $G_s$.
\end{proof}

Note that the reason that the blocks must be biconnected is to avoid edges that do not belong to cycles, and thus they have no role in the corresponding toric ideal.

\begin{rem1}{\rm
In graph theory, Theorem \ref{construction!} leads us to construct the family of graphs with $s\geq 2$ non bipartite blocks such that all its odd cycles share a common vertex. 

The only remaining open case is that the graph G has one non bipartite block such that $G$ has the strong odd cycle property. Here, the situation is completely different. A similar idea to the one applied above with the notion of a link vertex is the notion of an odd cycle transversal $D$. We recall that in graph theory, an odd cycle transversal of an undirected graph is a set of vertices of the graph that has a non empty intersection with every odd cycle in the graph. Removing the vertices of an odd cycle transversal from a graph leaves a bipartite graph. The problems that arise in our case is first the non-uniqueness of the odd cycle transversal and second the graph $G\setminus D$ is not always connected$;$ a remark that interrupted us in applying similar ideas. To the best of our knowledge, in graph theory, there are no algorithms that describe the family of graphs with the strong odd cycle property. 
}
\end{rem1}


\end{document}